\documentclass[final]{siamltex} 
\usepackage{amsmath} 
\usepackage{amsfonts} 
\usepackage{amssymb} 
\usepackage{amscd} 
\usepackage{mathrsfs} 
 \usepackage{graphicx} 
\usepackage{enumerate} 
  
\newtheorem{thm}{Theorem}[section] 
\newtheorem{cor}[thm]{Corollary} 
\newtheorem{lem}[thm]{Lemma} 
\newtheorem{prop}[thm]{Proposition}

\newtheorem{defn}[thm]{Definition}%[section] 

\newtheorem{example}[thm]{Example} 
\newtheorem{assumption}[thm]{Assumption} 

\newtheorem{remark}{Remark} 

\newcommand{\Z}{\mathbb{Z}} 
\newcommand{\R}{\mathbb{R}}

\newcommand{\GAC}{Global Attractor Conjecture~} 
\newcommand{\GACnospace}{Global Attractor Conjecture} 
 
\def\S{\mathcal S} 
\def\C{\mathcal C} 
\def\Re{\mathcal R} 
 
\title{The dynamics of weakly reversible population processes near facets\thanks{David Anderson was supported through grant 
NSF-DMS-0553687.  Anne Shiu was supported by a Lucent Technologies 
Bell Labs Graduate Research Fellowship.}}

\author{David F. Anderson\thanks 
{Department of Mathematics, University of Wisconsin at Madison, Madison, WI 53706 
({\tt anderson@math.wisc.edu}).}  
\and Anne Shiu\thanks 
{Department of Mathematics, University of California at Berkeley, Berkeley, CA 94720 
({\tt annejls@math.berkeley.edu}).}}

\begin{document} 
\slugger{siap}{2009}{}{}{} 
 
\maketitle  
 
%\bibliographystyle{plain} 
 
% : ABSTRACT 
\begin{abstract} 
  This paper concerns the dynamical behavior of weakly reversible, 
  deterministically modeled population processes near the facets 
  (codimension-one faces) of their invariant manifolds and proves that 
  the facets of such systems are ``repelling.'' It has been 
  conjectured that any population process whose network graph is 
  weakly reversible (has strongly connected components) is persistent. 
  We prove this conjecture to be true for the subclass of weakly 
  reversible systems for which only facets of the invariant manifold 
  are associated with semilocking sets, or siphons.  An important 
  application of this work pertains to chemical reaction systems that 
  are complex-balancing.  For these systems it is known that within 
  the interior of each invariant manifold there is a unique 
  equilibrium.  The Global Attractor Conjecture states that each of 
  these equilibria is globally asymptotically stable relative to the 
  interior of the invariant manifold in which it lies.  Our results 
  pertaining to weakly reversible systems imply that this conjecture 
  holds for all complex-balancing systems whose boundary equilibria 
  lie in the relative interior of the boundary facets.  As a 
  corollary, we show that the Global Attractor Conjecture holds for 
  those systems for which the associated invariant manifolds are 
  two-dimensional. 
 
  \vskip 0.1cm 
 
 
  \end{abstract} 
 
%: INTRODUCTION 
\section{Introduction} 
 
Population processes are mathematical models that describe the time 
evolution of the abundances of interacting ``species.''  To name a few 
examples, population processes can be used to describe the dynamics of 
animal populations, the spread of infections, and the evolution of 
chemical systems.  In these examples, the constituent species are the 
following: types of animals, infected and non-infected individuals, 
and chemical reactants and products, respectively.  How best to model 
the dynamics of a population process depends upon the abundances of 
the constituent species.  If the abundances are low, then the 
randomness of the interactions among the individual species is crucial 
to the system dynamics, so the process is most appropriately modeled 
stochastically.  % as a continuous-time Markov 
%chain.  
On the other hand, if the abundances are sufficiently high so that the 
randomness is averaged out at the scale of concentrations, then the 
dynamics of the concentrations can be modeled deterministically.  For 
precise statements regarding the relationship between the two models, 
see \cite{Kurtz1981, Kurtz72}.  In the present paper we consider 
deterministic models.  Also, we shall adopt the language associated 
with (bio)chemical reaction systems, which form a class of dynamical 
systems that arise in systems biology, and simply note that our 
results apply to any population process that satisfies our basic 
assumptions. 
 
The present work builds upon the body of work (usually called 
``chemical reaction network theory'') that focuses on the qualitative 
properties of chemical reaction systems and, in particular, those 
properties that are {\em independent of the values of the system 
  parameters}. See for example \cite{FeinDefZeroOne, Fein79, Guna, 
  HornJackson}.  Examples of reaction systems from biology include 
pharmacological models of drug interaction \cite{LRAT}, T-cell signal 
transduction models \cite{ChavezThesis,T-cell,Sontag01}, and enzymatic 
mechanisms \cite{SiegelMacLean}. This line of research is important 
because there are many % as there are potentially millions of 
biochemical reaction systems that may warrant study at one time or 
another, and these systems are typically complex and highly nonlinear. 
Further, the exact values of the system parameters are often unknown, 
and, worse still, these parameter values may vary from cell to cell. 
However, in a way that will be made precise in Section 
\ref{sec:model}, the network structure of a given system induces 
differential equations that govern its 
dynamics, %(up to parameter values), 
and it is this association between network structure and dynamics that 
can be utilized %to determine qualitative properties of 
% the dynamical systems  
without the need for detailed knowledge of parameter values. 
%maybe use the words ``systems biology'' in the above paragraph 
 
To introduce our main results, we recall three terms from the 
literature that will be defined more precisely later.  First, a 
directed graph is said to be {\em weakly reversible} if each of its 
connected components is strongly connected.  The directed graphs we 
consider in this paper are chemical reaction diagrams in which the 
arrows denote possible reactions and the nodes are linear combinations 
of the species which represent the sources and products of the 
reactions.  Second, for the systems in this paper a given trajectory 
is confined to an invariant polyhedron, which we shall denote by $P$. 
Such a polyhedron is called a {\em positive stoichiometric 
  compatibility class} in the chemical reaction network theory 
literature and the faces of its boundary are contained in the boundary 
of the positive orthant.  Third, {\em semilocking sets}, or {\em 
  siphons} in the Petri net literature \cite{Sontag07, 
  PetriNet,siphons}, are subsets of the set of species that 
characterize which faces of the boundary of $P$ allow for the 
existence of equilibria and $\omega$-limit points. 
 
The main result of this paper, Theorem \ref{thm:main}, concerns the 
dynamics of weakly reversible chemical reaction systems near {\em 
  facets} of $P$; a facet is a codimension-one face of $P$. 
% , in other words, its dimension is one less than that of $P$. 
Informally, Theorem~\ref{thm:main} states that weak reversibility of 
the reaction diagram guarantees the following: {\em for each point $z$ 
  found within the interior of a facet of $P$, there exists an open 
  (relative to $P$) neighborhood of $z$ within which trajectories are 
  forced away from the facet}.  Thus, Theorem \ref{thm:main} shows 
that weak reversibility guarantees that all facets are ``repelling.'' 
We will prove this theorem by demonstrating that for each facet there 
must exist a reaction that pushes the trajectory away from that facet 
and that the corresponding reaction rate dominates all others. 
 
The main qualitative results of this paper concern the long term 
behavior of systems, and as such we are interested in the set of 
$\omega$-limit points (accumulation points of trajectories).  A 
bounded trajectory of a dynamical system for which $\R^N_{\ge 0}$ is 
forward invariant is said to be {\em persistent} if no $\omega$-limit 
point lies on the boundary of the positive orthant.  Thus, persistence 
corresponds to a non-extinction requirement.  It has been conjectured 
that weak reversibility of a chemical reaction network implies that 
trajectories are persistent (for example, see \cite{FeinDefZeroOne}). 
Theorem \ref{thm:main} allows us to prove our main qualitative result, 
Theorem \ref{thm:persFacet}, which shows that this conjecture is true 
for the subclass of weakly reversible systems for which only facets of 
the invariant manifold are associated with semilocking sets.  We also 
point out in Corollary \ref{cor:dyn_non} that a slight variant of our 
proof of Theorem \ref{thm:main} shows that semilocking sets associated 
with facets are ``dynamically non-emptiable'' in the terminology of 
D.~Angeli {\em et al.} \cite{Sontag07}, thereby providing a large 
class of dynamically non-emptiable semilocking sets.

An important application of our main results pertains to chemical 
reaction systems that are {\em detailed-balancing} or, more generally, 
{\em complex-balancing} \cite{HornJackson}; these terms will be 
defined in Section~\ref{sec:GAC}.  For such systems, it is known that 
there is a unique equilibrium within the interior of each positive 
stoichiometric compatibility class $P$.  This equilibrium is called 
the Birch point in \cite{TDS} due to the connection to Birch's Theorem 
in Algebraic Statistics \cite[Section 2.1]{ASCB}.  Moreover, a strict 
Lyapunov function exists for this point, so local asymptotic stability 
relative to $P$ is guaranteed \cite{Fein79, HornJackson}.  An open 
question is whether all trajectories with initial condition in the 
interior of $P$ 
% (we call these {\em interior trajectories}) 
converge to the unique Birch point of $P$.  The assertion that the 
answer is `yes' is the content of the Global Attractor Conjecture 
%(GAC)  
\cite{TDS, HornJackson}.

The \GAC is a special case of the conjecture discussed earlier that 
pertains to weakly reversible systems; in other words, one must show 
that all complex-balancing systems are persistent 
\cite{FeinDefZeroOne}.  It is known that the set of $\omega$-limit 
points of such systems is contained within the set of equilibria 
\cite{ChavezThesis,Sontag01}, so the conjecture is equivalent to the 
statement that any equilibrium on the boundary of $P$ is not an 
$\omega$-limit point of an interior trajectory.  Recent work has shown 
that certain boundary equilibria are not $\omega$-limit points of 
interior trajectories.  For example, vertices of a positive 
stoichiometric compatibility class $P$ are not $\omega$-limit points 
of interior trajectories even if they are equilibria 
\cite{Anderson08,TDS}.  In addition, the \GAC recently has been shown 
to hold in the case that the system is detailed-balancing, $P$ is 
two-dimensional, and the system is conservative (meaning that $P$ is 
bounded) \cite{TDS}.  It is known that the underlying network of any 
detailed- or complex-balancing system necessarily is weakly 
reversible.  Therefore, all of the results of this paper apply in this 
setting and when combined with previous results \cite{Anderson08,TDS}, 
give our main contribution to the \GACnospace, Theorem 
\ref{thm:facetVtx}: {\em the \GAC holds for systems for which the 
  boundary equilibria are confined to facet-interior points and 
  vertices of $P$}.  %Thus, the present paper is complementary to both 
% \cite{Anderson08} and \cite{TDS}, which gave results concerned with 
% vertices of $P$ for complex-balancing systems. 
As a direct corollary to Theorem \ref{thm:facetVtx}, we can conclude 
that the \GAC holds for all systems for which $P$ is two-dimensional, 
in other words, a polygon, thereby extending the result in \cite{TDS}. 
 
We now describe the layout of the paper.  Section~\ref{sec:model} 
develops the mathematical model used throughout this paper.  In so 
doing we also present concepts from polyhedral geometry (Section 
\ref{sec:polyhedron}) that will be useful to us and formally define 
the notion of persistence (Section \ref{sec:persistence}).  In 
addition, the concept of a semilocking set is recalled.  Our main 
results are then stated and proven in Section~\ref{sec:MainThm}. 
Applications of this work to the \GAC is the topic of Section 
\ref{sec:GAC}.  Finally, Section~\ref{sec:exs} provides examples that 
illustrate our results within the context of related results. 
 
%:Section2 
\section{Mathematical formulation} 
\label{sec:model} 
 
In Sections \ref{sec:networks} and \ref{sec:kinetics}, we develop the 
mathematical model used in this paper and provide a brief introduction 
to chemical reaction network theory.  In Section~\ref{sec:polyhedron}, 
we present useful concepts from polyhedral geometry.  In 
Section~\ref{sec:persistence}, we recall the notions of persistence 
and semilocking sets.  Throughout the following sections, we adopt the 
notation $[n]:=\{ 1, 2, \dots, n\}$, for positive integers $n \in 
\Z_{>0}$. 
 
\subsection{Chemical reaction networks and basic terminology} 
\label{sec:networks} 
 
An example of a chemical reaction is denoted by the following: 
\begin{align*} 
  2X_{1}+X_{3} ~\rightarrow~ X_{2}~. 
\end{align*} 
The $X_{i}$ are called chemical {\em species} and $2X_{1}+X_{3}$ and 
$X_{2}$ are called chemical {\em complexes.}  Assigning the {\em 
  source} (or reactant) complex $2X_{1}+X_{3}$ to the vector $y = 
(2,0,3)$ and the {\em product} complex $X_{2}$ to the vector 
$y'=(0,1,0)$, we can write the reaction as $ y \rightarrow y' ~. $ In 
general we will denote by $N$ the number of species $X_{i}$, and we 
consider a set of $R$ reactions, each denoted by 
\begin{align*} 
  y_{k} \rightarrow y_{k}'~,  
\end{align*}  
for $k \in [R]$, and vectors $y_k, y_k' \in \Z^N_{\ge 0}$, with $y_k \ne 
y_k'$.  Note that if $y_k = \vec 0$ or $y_k' = \vec 0$, then this 
reaction represents an input or output to the system.  Note that any 
complex may appear as both a source complex and a product complex in 
the system.  For ease of notation, when there is no need for 
enumeration we typically will drop the subscript $k$ from the notation 
for the complexes and reactions. 
 
% For a given reaction, the $k$th say, there exists vectors $y_k, y_k' 
% \in \Z^N_{\ge 0}$, with $y_k \ne y_k'$, representing the number of 
% molecules of each species consumed and created, respectively, in one 
% instance of that reaction.  Using a slight abuse of notation, we 
% associate each such $y_k$ (and $y_k'$) with a linear combination of 
% the species in which the coefficient of $X_i$ is $y_{ki}$, the $i$th 
% component of $y_k$.% 
% Under this association, each $y_k$ (and $y_k'$) is termed a 
% {\em complex} of the system.  We may now denote any reaction by 
% the notation $y_k \to y_k'$, where $y_k$ is the source, or reactant, 
% complex and $y_k'$ is the product complex. 
\begin{defn} 
  Let $\S = \{X_i\}$, $\C = \{y\},$ and $\Re = \{y \to y'\}$ denote 
  sets of species, complexes, and reactions, respectively.  The triple 
  $\{\S, \C, \Re\}$ is called a {\em chemical reaction network}. 
  \label{def:crn} 
\end{defn}

To each reaction network, $\{\mathcal{S},\mathcal{C},\mathcal{R}\}$, 
we assign a unique directed graph (called a {\em reaction diagram}) 
constructed in the following manner.  The nodes of the graph are the 
complexes, $\mathcal{C}$.  A directed edge $(y,y')$ exists if and only 
if $y \to y'$ is a reaction in $\mathcal{R}$.  Each connected 
component of the resulting graph is termed a {\em linkage class} of 
the reaction diagram.  %, and we denote the number of linkage classes by~$l$. 
\begin{defn} 
  The chemical reaction network is said to be {\em weakly 
    reversible} if each linkage class of the corresponding reaction 
  diagram is strongly connected.  A network is said to be 
  {\em reversible} if $y' \to y \in \Re$ whenever $y \to y' \in 
  \Re.$ Later we will say that a chemical reaction system is 
  {\em weakly reversible} if its underlying network is. 
\end{defn} 
 
Let $x(t) \in \R^N$ denote the concentration vector of the species at 
time $t$ with initial condition $x(0)=x^0$.  We will show in 
Section~\ref{sec:kinetics} that the vector $x(t) - x^0$ remains within 
the span of the {\em reaction vectors} $\{y_k' - y_k\}$, i.e. in the 
linear space $S = \text{span}\{y_k' - y_k\}_{k \in [R]},$ for all 
time.  We therefore make the following definition. 
 
\begin{defn} 
  The {\em stoichiometric subspace} of a network is the linear 
  space $S = \text{span}\{y_k' - y_k\}_{k \in [R]}$. 
  \label{def:stoich_sub} 
\end{defn} 
 
It is known that under mild conditions on the rate functions of a 
system (see Section~\ref{sec:kinetics}), a trajectory $x(t)$ with 
strictly positive initial condition $x^0 \in \R^N_{>0}$ remains in the 
strictly positive orthant $\R^N_{>0}$ for all time (see Lemma~2.1 of 
\cite{Sontag01}).  Thus, the trajectory remains in the open set $(x^0 
+ S) \cap \mathbb{R}^N_{> 0}$, where $x^0 + S := \{z \in \R^N \ | \ z 
= x^0 + v, \text{ for some } v \in S\}$, for all time.  In other 
words, this set is {\em forward-invariant} with respect to the 
dynamics.  We shall refer to the closure of $(x^0 + S) \cap 
\mathbb{R}^N_{> 0}$, namely 
\begin{align} 
  P~:=~(x^0 + S) \cap \mathbb{R}^N_{\geq 0}~, 
  \label{posClass} 
\end{align} 
as a {\em positive stoichiometric compatibility class}.  We note that 
this notation is slightly nonstandard, as in previous literature it 
was the interior of $P$ that was termed the positive stoichiometric 
compatibility class.  In the next section, we will see that $P$ is a 
polyhedron. 
 
\begin{remark} \label{rmk:1} In spite of the notation, $P$ clearly 
  depends upon a choice of $x^0 \in \R^N_{>0}$.  Throughout the paper, 
  a reference to $P$ assumes the existence of a positive initial 
  condition $x^0 \in \R^N_{>0}$ for which $P$ is defined by 
  \eqref{posClass}. 
\end{remark} 
 
It will be convenient to view the set of species $\S$ as 
interchangeable with the set $[N]$, where $N$ denotes the number of 
species.  Therefore, a subset of the species, $W \subset \S$, is also 
a subset of $[N]$, and we will refer to the {\em $W$-coordinates} 
of a concentration vector $x \in \mathbb{R}^{N}$, meaning the 
concentrations $x_{i}$ for species $i$ in $W$.  Further, we will write 
$i \in W$ or $i \in [N]$ to represent $X_i \in W$ or $X_i \in \S$, 
respectively.  Similarly, we sometimes will consider subsets of the 
set of reactions $\mathcal{R}$ as subsets of the set $[R]$. 
 
\begin{defn} 
  The {\em zero-coordinates} of a vector $w \in \R^N$ are the 
  indices $i$ for which $w_i = 0$.  The {\em support} of $w$ is the 
  set of indices for which $w_i \ne 0$. 
  \label{def:support} 
\end{defn} 
 
Based upon Definition \ref{def:support} and the preceding remarks, 
both the set of zero-coordinates and the support of a vector $w$ can, 
and will, be viewed as subsets of the species.

% We will define the interior of $P$ and will review relevant concepts 
% of polyhedral geometry in Section~\ref{sec:polyhedron}.  This view 
% will justify naming $P$ an ``invariant polyhedron'' in \cite{TDS}. 
 
% Given that trajectories remain in their positive stoichiometric 
% compatibility classes $P$ for all time, we see that it is 
% appropriate to ask about the existence and stability of equilibria 
% of (\ref{eq:main}) within and relative to a positive compatibility 
% class.  We will take this viewpoint in Section~\ref{sec:GAC}. 
 
\subsection{The dynamics of a reaction system} 
\label{sec:kinetics} 
 
A chemical reaction network gives rise to a dynamical system by way of 
a rate function for each reaction. 
% In order to make a chemical reaction network a dynamical system we 
% must provide rate functions for each reaction. 
In other words, for each reaction $y_k \to y_k'$ we suppose the existence 
of a continuously differentiable function $\displaystyle R_k(\cdot) = 
R_{y_k \to y_k'}(\cdot)$ that satisfies the following assumption. 
\begin{assumption} For $k \in [R]$, $\displaystyle R_k(\cdot) = R_{y_k 
    \to y_k'}(\cdot): \R^N_{\ge 0} \to \R$ satisfies: 
  \begin{enumerate} 
  \item $R_{y_k \to y_k'}(\cdot)$ depends explicitly upon $x_i$ only 
    if $\displaystyle y_{ki} \ne 0$. 
  \item $\displaystyle \frac{\partial}{\partial x_i}R_{y_k \to 
      y_k'}(x) \ge 0$ for those $x_i$ for which $\displaystyle y_{ki} 
    \ne 0$, and equality can hold only if $x \in \partial \mathbb{R}^N_{\ge 
      0}$. 
  \item $R_{y_k \to y_k'}(x) = 0$ if $x_i = 0$ for some $i$ with 
    $\displaystyle y_{ki} \ne 0$. 
  \item If $\displaystyle 1 \le y_{ki} < y_{\ell i}$, then 
    $\displaystyle \lim_{x_i \to 0} \frac{R_{\ell}(x)}{R_{k}(x)} = 0$, 
    where all other $x_j > 0$ are held fixed in the limit. 
  \end{enumerate} 
  \label{assump:rates} 
\end{assumption} 
 
The final assumption simply states that if the $l$th reaction demands 
strictly more molecules of species $X_i$ as inputs than does the $k$th 
reaction, then the rate of the $l$th reaction decreases to zero faster 
than the $k$th reaction, as $x_i \to 0$. The functions $R_{k}$ are 
typically referred to as the {\em kinetics} of the system and the 
dynamics of the system are given by the following coupled set of 
nonlinear ordinary differential equations: 
 
\begin{equation} 
  \dot x(t) = \sum_{k \in [R]} R_{k}(x(t))(y_k' - y_k)~. 
  \label{eq:main_general} 
\end{equation} 
Integrating \eqref{eq:main_general} yields 
\begin{equation*} 
  x(t) = x^0 + \sum_{k \in [R]} \left(\int_0^t R_k(x(s)) ds \right) 
  (y_k' - y_k)~.  
\end{equation*} 
Therefore, $x(t) - x^0$ remains in the stoichiometric subspace, $S = 
\text{span}\{y_k' - y_k\}_{k \in [R]}$, for all time, confirming the 
assertion made in the previous section. 
 
The most common kinetics, and the choice we shall make throughout the 
remainder of this paper, is that of {\em mass action kinetics}. A 
chemical reaction system is said to have mass action kinetics if all 
functions $R_{k}$ take the following multiplicative form: 
\begin{equation} 
  R_{k}(x) =  \kappa_k x_1^{y_{k1}} x_2^{y_{k2}} \cdots x_N^{y_{kN}} 
  =: \kappa_k x^{y_k}~, 
  \label{eq:massaction} 
\end{equation} 
for some positive reaction rate constants $\kappa_{k}$, where we have 
adopted the convention that $0^0 = 1$ and the final equality is a 
definition.  It is easily verified that each $R_k$ defined via 
\eqref{eq:massaction} satisfies Assumption \ref{assump:rates}. 
Combining \eqref{eq:main_general} and \eqref{eq:massaction} gives the 
following system of differential equations: 
\begin{equation} 
  \dot x(t) = \sum_{k \in [R]} \kappa_k x(t)^{y_k}(y_k' - y_k) =: f(x(t))~, 
  \label{eq:main} 
\end{equation} 
where the last equality is a definition.  This dynamical system is the 
main object of study in this paper. 
 
A concentration vector $\overline x \in \R^N_{\geq 0}$ is an {\em 
  equilibrium} of the mass action system \eqref{eq:main} if 
$f(\overline x) = 0$.  Given that trajectories remain in their 
positive stoichiometric compatibility classes $P$ for all positive 
time, we see that it is appropriate to ask about the existence and 
stability of equilibria of system (\ref{eq:main}) within and relative 
to a positive stoichiometric compatibility class $P$.  We will take 
this viewpoint in Section~\ref{sec:GAC}. 
 
\begin{remark} 
  We note that every result in this paper holds for any chemical 
  reaction systems with kinetics that satisfy Assumption 
  \ref{assump:rates}.  Nonetheless, we choose to perform our analysis 
  in the mass action case for clarity of exposition. 
\end{remark} 
 
\subsection{Connection to polyhedral geometry} 
\label{sec:polyhedron} 
 
We now recall terminology from polyhedral geometry that will be 
useful; we refer the reader to the text of G. Ziegler for further 
details \cite{Ziegler}. 
\begin{defn} 
  The {\em half-space} in $\mathbb{R}^m$ defined by a vector $v \in 
  \mathbb{R}^m$ and a constant $c \in \mathbb{R}$ is the set 
  \begin{align} 
    H_{v,c} ~:=~ \left\{ x \in \mathbb{R}^m ~|~ \langle v, x \rangle 
      \geq c \right\}~. 
  \end{align} 
  % We see that $H_{v,c}$ is the set of all points lying on one side 
  % of the hyperplane defined by $\langle v, - \rangle = c$. 
  A (convex) {\em polyhedron} in $\mathbb{R}^m$ is an intersection of 
  finitely many half-spaces. 
\end{defn} 
 
For example, the non-negative orthant $\mathbb{R}^N_{\geq 0}$ is a 
polyhedron, as it can be written as the intersection of the $N$ 
half-spaces $H_{e_i,0}$, where the $e_i$'s are the canonical unit 
vectors of $\mathbb{R}^N$. We now give three elementary facts about 
polyhedra from which we will deduce the fact that positive 
stoichiometric compatibility classes $P$ are polyhedra.  First, any 
linear space of $\mathbb{R}^m$ is a polyhedron.  Second, any 
translation $x + Q$ of a polyhedron $Q$ by a vector $x \in 
\mathbb{R}^m$ is again a polyhedron.  Third, the intersection of two 
polyhedra is a polyhedron.  Therefore, as a translate $(x^0 + S)$ and 
the orthant $\mathbb{R}^N_{\geq 0}$ are both polyhedra, it follows 
that the positive stoichiometric compatibility class $P$ defined by 
(\ref{posClass}) is indeed a polyhedron. 
 
We continue with further definitions, which will allow us later to 
discuss {\em boundary equilibria} (those equilibria of (\ref{eq:main}) 
on the boundary of 
$P$).  %It can be shown that the zero-coordinate set of any boundary 
       %equilibrium is a semilocking set (see next Section~- this term 
       %is not yet defined). 
\begin{defn} 
  Let $Q$ be a polyhedron in $\mathbb{R}^m$.  The {\em interior} of 
  $Q$, denoted by $\operatorname{int}(Q)$, is the largest relatively open subset 
  of $Q$.  The {\em dimension} of $Q$, denoted by $\dim (Q)$, 
  is %the dimension of its affine span; concretely, it is 
  the dimension of the span of the translate of $Q$ that contains the 
  origin. % of $\mathbb{R}^m$. 
\end{defn} 
 
For example, the dimension of $P$ equals the dimension of the 
stoichiometric subspace $S$: $\dim(P)=\dim(S)$.  % b/c c^{0} is pos. vector 
We now define the 
faces of a polyhedron. 
%[page 3 of \cite{Ziegler}] 
\begin{defn} \label{def:face}  
Let $Q$ be a polyhedron in 
  $\mathbb{R}^m$.  For a vector $v\in \mathbb{R}^m$, the {\em face} of 
  $Q$ that it defines is the (possibly empty) set of points of $Q$ 
  that minimize the linear functional $\langle v, \cdot \rangle: 
  \mathbb{R}^m \rightarrow \mathbb{R}$. 
\end{defn}   
 
If the minimum in Definition \ref{def:face} (denoted $c_{\min}$) is 
attained, then we can write the face as $F=Q \cap H_{v,c_{\min}} \cap 
H_{-v,c_{\min}}$.  Therefore any face is itself a polyhedron, so we 
may speak of its dimension or its 
interior.  % $\operatorname{int}(F)$. 
\begin{defn}\label{def:facetVtx} 
Let $Q$ be a polyhedron in 
  $\mathbb{R}^m$.  A {\em facet} of $Q$ is a 
  face whose dimension is one less than that of $Q$.  A {\em vertex} 
  is a nonempty zero-dimensional face (thus, it is a point). 
\end{defn} 
 
We make some remarks.  First, note that what we call the ``interior'' 
is sometimes defined as the ``relative interior'' \cite{Ziegler}. 
Second, vertices are called ``extreme points'' in \cite{Anderson08}. 
Third, the interior of a vertex is seen to be the vertex itself. 
Fourth, the boundary of $Q$ is the disjoint union of the interiors of 
the proper faces of $Q$.  %By proper we mean... 
 
We now return to the positive stoichiometric classes $P$  \eqref{posClass}.  For a subset of the set of 
species $W \subset \S$, let $Z_{W} \subset \R^N$ denote its 
{\em zero set}: 
\begin{align*} 
  Z_W = \{x \in \R^{N} ~:~ x_i = 0 \text{ if } i \in W\}~. 
\end{align*} 
It can be seen that for any face $F$ of a positive stoichiometric 
class $P$, there exists some possibly non-unique subset  
$W \subset \S$ such that 
\begin{align}\label{faceOfP} 
  F = F_W := P \cap Z_W~. 
\end{align} 
In other words, each face of $P$ is the set of points of $P$ whose set 
of zero-coordinates contains a certain subset $W \subset \S$.  However, 
it is important to note that for some subsets $W$, the face is empty: 
$F_W = \emptyset$, and therefore 
% if $P$ does not intersect the face of the positive orthant 
% associated with $W$, and 
no nonempty face of $P$ corresponds with such a $W$.  In this case we 
say that the set $Z_W$ is {\em stoichiometrically unattainable}. 
We see also that $F_{W}=P$ if and only if $W$ is empty. For 
definiteness, if there exist subsets $W_1 \subsetneqq W_2 \subset \S$  
%with $W_1 \ne W_2$  
for which $F_{W_1} = F_{W_2}$, we denote the face 
by $F_{W_2}$. 
Under this convention, it can be seen that the interior of a face $F_W$ is 
\begin{align} 
  \label{faceInt} 
  \operatorname{int}(F_W)=\left\{ ~x\in P ~|~ x_i=0 \text{ if and only 
      if } i \in W ~ \right\} . 
\end{align} 
We remark that the set $\text{int}(F_W)$ was denoted by $L_W \cap P$ 
in \cite{Anderson08}. 
 
The following example illustrates the above concepts. We note that in 
the interest of clarity we denote species by $A, B, C, \dots$ 
rather than $X_{1}, X_{2}, X_{3}, \dots$ in all examples. 
 
% :Example of polyhedral concepts via P 
\begin{example} Consider the chemical reaction system which arises from the following reaction 
  diagram: 
  \begin{align}\label{reactionEx1} 
    2A ~ \underset{\kappa_2}{\overset{\kappa_1}{\rightleftarrows}} ~ 
    A+B \quad , \quad B ~ 
    \underset{\kappa_4}{\overset{\kappa_3}{\rightleftarrows}} ~ C ~, 
  \end{align} 
  where we use the standard notation of labeling a reaction arrow by 
  the corresponding reaction rate constant.  The stoichiometric 
  subspace $S$ in $\mathbb{R}^3$ is spanned by the two reaction 
  vectors $(-1,1,0)$ and $(0,-1,1)$.  A positive stoichiometric 
  compatibility class is depicted in Figure \ref{Fig:ex1}; it is a 
  two--dimensional {\em simplex} (convex hull of three affinely 
  independent points, in other words, a triangle) given by 
  \begin{align} 
    P~=~ \left\{~ (x_a, x_b, x_c ) \in \mathbb{R}^3_{\geq 0} ~|~ 
      x_a+x_b + x_c = T ~ \right\}~, 
    \label{eq:triangle} 
  \end{align} 
  for positive total concentration $T>0$. 
%Figure 
% 
  \begin{figure}%[htbp] 
    \begin{center} 
      \includegraphics[scale=0.40]{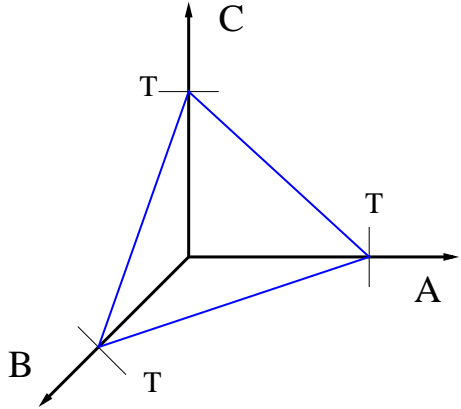}%{figure.eps}%plane2.pdf} 
      \caption{Positive stoichiometric compatibility class $P$ for 
        chemical reaction system (\ref{reactionEx1}).} 
      \label{Fig:ex1} 
    \end{center} 
  \end{figure} 
  The three facets (edges) of each positive stoichiometric 
  compatibility class $P$ are one-dimensional line segments: 
  \begin{align*} 
    F_{\{A\}} &= \left\{~ (0, x_b, x_c ) \in \mathbb{R}^3_{\geq 0} ~|~ 
      x_b + x_c = T ~ \right\}~,\\ 
    F_{\{B\}} &= \left\{~ (x_a, 0, x_c ) \in \mathbb{R}^3_{\geq 0} ~|~ 
      x_a + x_c = T ~ \right\}~,\\ 
    F_{\{C\}}&= \left\{~ (x_a, x_b, 0 ) \in \mathbb{R}^3_{\geq 0} ~|~ 
      x_a + x_b = T ~ \right\}~, 
  \end{align*} 
  and the three vertices are the three points $F_{ \{ A, B \} } = 
  \{(0,0,T)\}$, $F_{ \{ A, C \} } = \{ (0,T,0) \} $, and $F_{ \{ B, C 
    \} } = \{ (T,0,0) \} $.  Finally, the set $Z_{\{A,B,C\}} = 
  \{(0,0,0)\}$ is stoichiometrically unattainable.  We will revisit 
  this reaction network in Example \ref{ex:1}. 
  \label{ex:P} 
\end{example} 
 
\subsection{Persistence and semilocking sets} 
\label{sec:persistence} 
 
% By Lemma~\ref{lemma:pos}, each trajectory must remain within 
% $\R^m_{>0}$ if its initial condition is in $\R^m_{>0}$; therefore 
% the linear subsets of interest are the intersections of the 
% stoichiometric compatibility classes and $\R^m_{>0}$.  Recall that 
% in the introduction these sets were termed the positive 
% stoichiometric compatibility classes. 

% :Define omega-limit points 
Let $x(t)$ be a solution to~\eqref{eq:main} with strictly positive 
initial condition $x^0 \in \R^N_{>0}$.  The set of 
{\em $\omega$-limit points} for this trajectory is the set of 
accumulation points: 
\begin{equation} 
  \omega(x^0)  ~:=~  \{~x \in \R^N_{\ge 0} ~|~ x(t_n) \to {x,}  \text{ 
    for some sequence } t_n \to \infty ~ \}. %subset of P 
  \label{omega} 
\end{equation} 
% important - we define omega-limit points only for interior 
% trajectories. 
 
\begin{defn} 
  A bounded trajectory with initial condition $x^0$ is said to be {\em 
    persistent} if $\omega(x^{0}) \cap \partial \R^{N}_{\ge 0} = 
  \emptyset$.  A dynamical system with bounded trajectories is {\em 
    persistent} if each trajectory with strictly positive initial 
  condition is persistent. 
  \label{def:persistence} 
\end{defn} 
 
In order to show that a system is persistent, we must understand which 
points on the boundary of a positive stoichiometric class are capable 
of being $\omega$-limit points.  To this end, we recall the following 
definition from the literature. 
 
\begin{defn} 
  A nonempty subset $W$ of the set of species is called a\/ {\em 
    semilocking set} if for each reaction in which there is an element 
  of $W$ in the product complex, there is an element of $W$ in the 
  reactant complex. 
  % $W$ is called a\/ {\rm locking set} if every reactant complex 
  % contains an element of~$W$. 
  \label{def:locking} 
\end{defn} 
 
\begin{remark} 
  The notion of a semilocking set is the same as a {\em siphon} in 
  the Petri net literature.  See \cite{Sontag07,PetriNet, siphons}. 
\end{remark} 
 
The intuition behind semilocking sets lies in the following 
proposition, which is the content of Proposition 5.5 of D. Angeli {\em et al.}  
\cite{PetriNet}.  Related results that concern the ``reachability'' of species include  
Theorems 1 and 2 in the textbook of A. Vol\'{}pert and S. Khud\t{ia}ev \cite[Section 12.2.3]{Volp}. 
 
\begin{prop}[\cite{PetriNet}] 
  Let $W \subset \S$ be non-empty.  Then $W$ is a semilocking set if 
  and only if the face $F_W$ is forward invariant for the dynamics  
  \eqref{eq:main}. 
  \label{prop:semilocking_intution} 
\end{prop} 
 
The above result holds because semilocking sets are characterized by 
the following property: if no species of $W$ are present at time zero, 
then no species of $W$ can be produced at any time in the future.  In 
other words, these species are ``locked'' at zero for all time.  If in 
addition the reaction network is weakly reversible, then it is 
straightforward to conclude the following: if a linkage class has a 
complex whose support contains an element of $W$, then the rates of 
all reactions within that linkage class will be zero for all positive 
time: $\kappa_k x^{y_k}=0$ for any such reaction $y_k \to y_k' $.  In other words, certain linkage classes are ``shut off.'' 
 
In light of the characterization of the interior of a face $F_{W}$ 
given in (\ref{faceInt}), the following theorem is proven in 
\cite{Anderson08,PetriNet}; it states that the semilocking sets are 
the possible sets of zero-coordinates of boundary $\omega$-limit 
points. 
 
\begin{thm}[\cite{Anderson08,PetriNet}] 
  Let $W \subset \S$ be a nonempty subset of the set of species.  Let $ x^0 \in 
  \R^N_{>0}$ be a strictly positive initial condition for the  
  system~\eqref{eq:main}, and let $P=(x^0+S) \cap \mathbb{R}^N_{\geq 0}$ 
  denote the corresponding positive stoichiometric compatibility 
  class.  If there exists an $\omega(x^0)$-limit point, $z \in 
  \omega(x^0)$, and a subset of the species, $W$, such that $z$ is 
  contained within the interior of the face $F_W$ of $P$, then $W$ is 
  a semilocking set. 
  \label{thm:omegapoints} 
\end{thm} 
 
Theorem \ref{thm:omegapoints} will be used in conjunction with results in the next 
section to prove the persistence of the following class of weakly 
reversible systems: those for which each semilocking set $W$ 
satisfies $\text{dim}(F_W) = \text{dim}(P) - 1$ (and so $F_W$ is a 
facet of $P$) or $F_W = \emptyset$ (and so $Z_W$ is stoichiometrically 
unattainable); see Theorem~\ref{thm:persFacet}. 
 
% :Section of main results 
\section{Main results} 
\label{sec:MainThm} 
 
In order to state Theorem \ref{thm:main}, we need the following definition. 
 
%: REVISED definition of repelling 
\begin{defn} 
  Let $Q \subset P$ be an open set relative to $P$, for which 
  $\emptyset \ne Q \cap \partial P \subset F_W$, for some face $F_W$ of $P$. 
  Then the face   
  $F_W$ is {\em repelling in the neighborhood} $Q \cap \operatorname{int}(P)$ 
    with respect to the dynamics \eqref{eq:main} if 
  \begin{equation} 
         \sum_{i \in W} x_i f_i(x) ~\ge~ 0 
\label{eq:def_repelling} 
%  \frac{d}{dt} \text{dist} \left( x(t), F_{W} \right) \geq 0, 
  \end{equation} 
  for all $x \in Q \cap \operatorname{int}(P)$, where the $f_i$ are the functions given in (\ref{eq:main}). 
    \label{def:repelling} 
\end{defn} 
 
\begin{remark} 
  Note that $F_W$ is repelling in the neighborhood $Q \cap 
  \operatorname{int}(P)$ with respect to \eqref{eq:main} if and only 
  if $\frac{d}{dt} \operatorname{dist} \left( x(t), F_{W} \right) \geq 
  0$ whenever $x(t) \in Q\cap \operatorname{int}(P)$.  Thus, $F_W$ is 
  repelling in a neighborhood if any trajectory in the neighborhood 
  can not get closer to the face $F_W$ while remaining in the 
  neighborhood. 
\end{remark} 
 
The following theorem is our main technical result. 
%a neighborhood in which the facet is repelling.  
 
\begin{thm} 
  Let $\{\S, \C, \Re\}$ be a weakly reversible chemical reaction 
  network with dynamics governed by mass action kinetics 
  \eqref{eq:main}.  Let $W \subset \S$ be such that $F_W$ is a facet 
  of $P$, and take $z$ to be in the interior of $F_W$.  Then there 
  exists a $\delta > 0$ for which the facet $F_W$ is repelling in the 
  neighborhood $B_{\delta}(z) \cap \operatorname{int}(P)$, where 
  $B_{\delta}(z)$ is the open ball of radius $\delta$ centered at the 
  point $z$. 
  \label{thm:main} 
\end{thm} 
 
% :proof of main theorem 
\begin{proof} 
  For the time being we will assume that there is only one linkage 
  class in the reaction diagram.  The proof of the more general case 
  is similar and will be discussed at the end. Also we assume that 
  $W=\{X_1,\dots, X_M\}$ for some $M \le N$. 
    
  Letting $s:=\dim (S) = \dim(P)$, the facet $F_{W}$ has dimension 
  $s-1$, which, by definition, means that $Z_{W} \cap S$ is an 
  $(s-1)$-dimensional subspace of $S$.  Let $\pi ~:~ \mathbb{R}^{N} 
  \rightarrow \mathbb{R}^{N}$ be the projection onto the first $M$ 
  coordinates, that is, it is given by $\pi(x_1,x_2,\dots,x_M, x_{M+1},\dots,x_N) := 
  (x_1, x_2, \dots,x_M,0,\dots,0)$.  As a shorthand we will also write 
  $x|_W$ for $\pi(x)$.  Because $Z_W \cap S=\ker(\pi|_S)$ has dimension $(s-1)$, it 
  follows that the image $\pi(S)$ is one-dimensional. 
  Therefore, we may let $v \in S$ be such that $v|_{W}$ spans the 
  projection $\pi(S)$.  We also let $\{w_2,\dots,w_s\}$ span the subspace $Z_{W} 
  \cap S$ so that $\{v,w_2,\dots,w_s\}$ is a basis for the subspace $S$.  We simply 
  note for future reference that by construction we have 
  \begin{equation} 
    y'| _{W} - y|_{W} ~\in ~ \text{span}(v|_{W})~, 
    \label{diffVec} 
  \end{equation} 
  for each reaction $y \to y' \in \Re$.  Finally, for $x \in \R^N$ we 
  define $x|_{W^c}$ similarly to $x|_W$; that is, $x|_{W^c}$ is the 
  projection of $x$ onto the final $N-M$ coordinates. % defined by setting the  first $M$ coordinates to zero. 
 
  We may assume that all coordinates of $v|_{W}$ are non-zero, for 
  otherwise the concentrations of certain species $X_{j} \in W$ would remain 
  unchanged under the action of each reaction (note that we 
  necessarily have $w_{k \, i} = 0$ for $k \in \{2,\dots,s\}$ and $i 
  \in \{1,\dots,M\}$).  In such a case, the concentrations of those species $X_{j}$ would remain constant in time, so we could simply remove them from the system by incorporating 
  their concentrations into the rate constants appropriately. 
   
  We will show that $v|_{W}$ has coordinates all of one sign and will use 
  this fact to guarantee the existence of a ``minimal complex'' (with 
  respect to the elements of $W$).  We will then show that this 
  minimal complex corresponds with a dominating monomial that appears 
  as a positive term in each of the first $M$ components of 
  \eqref{eq:main}. %, and the result will be shown. 
   
  Suppose, in order to find a contradiction, that $v|_{W}$ has 
  coordinates of both positive and negative sign; that is, assume that 
  $v_{i} < 0 < v_{j}$ for some indices $i, \ j \le M$.  Let $u:= 
  v_{j}e_{i}-v_{i}e_{j} \in \mathbb{R}^{N}_{\ge 0}$ (where $e_{l}$ 
  denotes the $l$th canonical basis vector).  It follows that $u \in 
  S^{\perp}$ because $(i)$ $\langle u, v\rangle = 0$ by construction, 
  and $(ii)$ $\langle u, w_k\rangle = 0$ for all $k\in \{2,\dots,s\}$ 
  because these vectors have non-overlapping support.  Note also that 
  $\langle u, z \rangle = 0$ because the support of $u$ is a subset of $W$ whereas the support of $z$ is $W^c$.  Let $x^0 \in 
  \R^N_{>0} \cap P$ (such a point always exists by Remark~\ref{rmk:1}).   
  As $z$ and $x^{0}$ both lie in $P$, there exist 
  constants $\alpha_k \in \R$ for $k \in [s]$, such that 
  \begin{equation*} 
    z = x^0 + \alpha_1 v + \sum_{k =2}^s \alpha_k w_k~. 
  \end{equation*} 
  Combining all of the above we conclude that  
  \begin{equation*} 
    0  ~=~ \langle u, z \rangle ~=~ \langle u, x^0 \rangle + \alpha_1 \langle u, v\rangle + \sum_{k = 2}^s \alpha_k \langle u, w_k \rangle  
    = \langle u, x^{0} \rangle  
    ~>~ 0~, 
  \end{equation*} 
  where the final inequality holds because $u$ is non-negative and 
  nonzero and $x^0$ has strictly positive components.  This is a 
  contradiction, so we conclude that $v|_{W}$ does not have both 
  positive and negative coordinates, and, without loss of generality, 
  we now assume that all coordinates of $v|_{W}$ are positive. 
     
  We recall from \eqref{diffVec} that $y'| _{W} - y|_{W} ~\in ~ 
  \text{span}(v|_{W})$ for each reaction $y \to y' \in \Re$.  For 
  concreteness we let ${y'_k}|_W - {y_k}|_W = \gamma_k v|_{W}$ for 
  some $\gamma_k \in \R$ where $k \in [R]$.  Combining this with the 
  fact that $v_i > 0$ for each $i \in \{1,\dots,M\}$ shows that each 
  reaction yields either $(i)$ a net gain of all species of $W$, 
  $(ii)$ a net loss of all species of $W$, or $(iii)$ no change in any 
  species of $W$ and, moreover, that there exists a $\tilde{y} \in \C$ 
  such that $\tilde{y}|_W \preceq y|_W$ for all $y \in \C$, where we 
  say $x \preceq y$ for $x,y \in \R^N$ if $x_i \le y_i$ for each $i$. 
  Note that it is the sign of $\gamma_k$ that determines whether a 
  given reaction accounts for an increase or decrease in the 
  abundances of the elements of $W$.

  % We now write 
%   \begin{align*} 
%   x^{y}~=~ x|_{W}^{y|_{W}} \cdot x|_{W^c} ^{y|_{W^c}} 
%  \end{align*} 
%  for the usual monomial that corresponds to the complex $y$; this is 
%  the monomial that appears in the differential equations 
%  (\ref{eq:main}).  Note that our choice of notation separates the 
%  contribution of the $W$-species concentrations, which will be small 
%  near $z$, from that of the other species, which will be bounded 
%  both above and below near $z$. 
 
  We now find a neighborhood of positive radius $\delta$ around $z$, 
  $B_{\delta}(z)$, for which the facet $F_W$ is repelling in the 
  neighborhood $B_{\delta}(z) \cap \text{int}(P)$. The first condition 
  we impose on $\delta$ is that it must be less than the distance 
  between $z$ and any proper face of $P$ that is not $F_{W}$, which 
  can be done because $z$ is in the interior of the facet.  Also, this 
  condition ensures that for any point $x \in B_{\delta}(z) \cap 
  \text{int}(P)$, the coordinates $x_i$, for $i > M$, are uniformly 
  bounded both above and below.  Therefore, there exist constants 
  $D_{\text{min}}$ and $D_{\text{max}}$ such that for all $x \in 
  B_{\delta}(z) \cap \text{int}(P)$ and all complexes $y$, we have the 
  inequalities 
  \begin{align} 
    \label{D_bdd} 
    0~<~ D_{\text{min}} ~<~  x|_{W^c} ^{y|_{W^c}}~<~ 
    D_{\text{max}}~. 
  \end{align} 
   
  The monomial $x|_W^{\tilde{y}|_W}$ will dominate all other monomials 
  for $x \in B_{\delta}(z)\cap \text{int}(P)$ for sufficiently small 
  $\delta$, and this will force trajectories away from the facet. To 
  make this idea precise, let $R_{+}$ denote those reactions that 
  result in a net gain of the species in $W$ and $R_{-}$ those that 
  result in a net loss.  We will write $y_k \to y_k' \in R_{+}$ and $y_k 
  \to y_k' \in R_{-}$ to enumerate over those reactions.  We now have 
  that for $i \in [M]$ and $x \in B_{\delta}(z)\cap \text{int}(P)$ and 
  for sufficiently small $\delta>0$,  
  \begin{align} 
  \begin{split} 
    f_i(x) ~&=~ v_i \sum_{y_k \to y_k' \in R_{+}} \gamma_k \kappa_k x|_W^{{y_k}|_W}x|_{W^c}^{{y_k}|_{W^c}} ~-~ v_i \sum_{y_k \to y_k' \in R_{-}} \left | \gamma_k \right |\kappa_k x|_W^{{y_k}|_W}x|_{W^c}^{{y_k}|_{W^c}} \\ 
    ~& \ge~ v_i D_{\text{min}} \sum_{y_k \to y_k' \in R_{+}}\gamma_k 
    \kappa_k x|_W^{{y_k}|_W} ~-~ v_i D_{\text{max}} \sum_{y_k \to y_k' \in 
      R_{-}} \left |\gamma_k \right | \kappa_k 
    x|_W^{{y_k}|_W}~. 
    \end{split} 
     \label{eq:bigbound} 
  \end{align} 
  Finally, by weak reversibility (and possibly after choosing a 
  different $\tilde y$ that still satisfies the minimality condition), 
  there is a reaction $\tilde{y} \to y' \in \Re$ with  
  $\tilde{y}_i < y'_i$  
  %$\tilde{y}|_{W\,i} < y'|_{W\, i}$  
    for all $i \in \{1, \dots, M\}$ (i.e. the 
  $\gamma_k$ associated with this reaction is strictly positive). 
  This reaction, which belongs to $R_{+}$, has a monomial, 
  $x|_W^{\tilde{y}|_{W}}$, that necessarily dominates all monomials 
  associated with reactions in $R_{-}$ (which necessarily have 
  source complexes that contain a {\em higher} number of each 
  element of $W$ than $\tilde y$ does).  Combining this fact with 
  \eqref{eq:bigbound} shows that $f_i(x) \ge 0$ for $i \in [M]$ and for $x 
  \in B_{\delta}(z) \cap \text{int}(P)$,  
and therefore, the facet $F_W$ is repelling in the neighborhood $B_{\delta}(z) \cap \text{int}(P)$. 
   
  In the case of more than one linkage class, each linkage class will have its own minimal complex that 
  will dominate all other monomials associated with that linkage 
  class.  Thus the desired result follows. 
\end{proof} 
 
\begin{remark} 
  Note that weak reversibility was used in the previous proof only to 
  guarantee the existence of the reaction $\tilde{y} \to y'$, where 
  $\tilde{y}|_{W}$ is minimal and  
    $\tilde{y}_i < y'_i$  
  %$\tilde{y}|_{W\, i} < y'|_{W\, i}$ 
  for all $i \in \{1, \dots, M\}$, and was {\em not} needed to 
  prove the existence of such a $\tilde{y}$.  If the network were not 
  weakly reversible, but such a reaction nevertheless existed, then the same 
  proof goes through unchanged. 
\end{remark}

\begin{cor} 
  Let $\{\S, \C, \Re\}$ be a weakly reversible chemical reaction 
  network with dynamics governed by mass action kinetics 
  \eqref{eq:main} such that all trajectories are bounded. Suppose there exists a subset $W \subset \S$, a positive initial condition $x^0 \in \R^N_{>0}$, and a point $z 
  \in \omega(x^0) \cap F_W$ such that $F_W$ is a facet of $P$.  Then 
  $\omega(x^0) \cap \partial F_W \ne \emptyset$. 
  \label{cor:edge} 
\end{cor} 
 
\begin{proof} 
  Suppose not.  That is, suppose that $\omega(x^0)\cap F_W \subset 
  \text{int}(F_W)$ holds.  Let $\mathcal{Y} :=\omega(x^{0}) \cap F_{W} $. 
  We claim that $\mathcal{Y}$ is a compact set; indeed, the trajectory 
  $x(t)$ is bounded so $\mathcal{Y}$ is as well, and $\mathcal{Y}$ is 
  the intersection of two closed sets, and therefore is itself closed. 
   
  Combining the compactness of $\mathcal{Y} \subset \text{int}(F_W)$ 
  with Theorem~\ref{thm:main}, we obtain that there exists an open covering 
  of $\mathcal{Y}$ consisting of a finite number of balls 
  $B_{\delta_{i}}(z_{i})$ of positive radius $\delta_{i}$, each 
  centered around an element $z_{i}$ of $\mathcal{Y}$, such that $(i)$ 
  each $\delta_{i}$ is sufficiently small so that 
  $\overline{B_{\delta_i}(z_i)} \, \cap \, \partial F_W = \emptyset$, 
  and $(ii)$, for $Q := \cup_i B_{\delta_i}(z_i)$, the facet $F_W$ is repelling in 
  $Q \cap 
  \text{int}(P)$.  Combining these 
  facts with the existence of $z \in \omega(x^0) \cap \text{int}(F_W) 
  \cap Q$ shows the existence of a point $w \in \omega(x^0) \cap 
  \text{int}(F_W) \cap \partial Q$.  However, this is impossible because $w \in 
  \omega(x^0) \cap \text{int}(F_W) = \mathcal{Y}$ necessitates that $w 
  \in \mathcal{Y} \subset \text{int}(Q)$.  Thus, the claim is shown. 
\end{proof} 
 
We may now present our main qualitative result. 
 
\begin{thm} 
  Let $\{\S, \C, \Re\}$ be a weakly reversible chemical reaction 
  network with dynamics governed by mass action kinetics 
  \eqref{eq:main} such that all trajectories are bounded.  Suppose 
  that for each semilocking set $W$, the corresponding face $F_{W}$ is 
  either a facet ($\text{dim}(F_W) = \text{dim}(P) - 1$) or is empty. 
  Then the system is persistent. 
  \label{thm:persFacet} 
\end{thm} 
 
\begin{proof} 
  This is an immediate consequence of Theorem \ref{thm:omegapoints} and Corollary \ref{cor:edge}. 
\end{proof} 
 
\subsection{Connection to dynamically non-emptiable semilocking sets} 
\label{sec:dyn_non_emp} 
 
In \cite{Sontag07} the notion of a ``dynamically non-emptiable'' 
semilocking set is introduced.  If for a given semilocking set $W$ we 
define two sets  
\begin{align*} 
  C(W) & := \left\{ ~ 0 \preceq \alpha \in \R^R_{\ge 0} \, : \, w = \sum_{k=1}^R \alpha_k (y'_k - y_k) \text{ satisfies } w|_{W} \preceq 0~ \right\} \text{~and} \\ 
  \mathcal{F}_{\epsilon}(W) &:= \left\{~ 0 \preceq \alpha \in \R^R_{\ge 0} 
  \, : \, \alpha_j \le \epsilon \alpha_i,~ \forall i,j \in [R] \text{ 
    such that } {y_i}|_{W}  \precneqq {y_j}|_{W}  ~ \right\}~, 
\end{align*} 
then $W$ is said to be {\em dynamically non-emptiable} if $C(W) \, 
\cap \, \mathcal{F}_{\epsilon}(W) = \{0\}$ for some $\epsilon > 0$.  Here the notation $z\preceq z'$ means that all coordinates satisfy the inequality $z_i \leq z_i'$ and $z \precneqq z'$ means that furthermore at least one inequality is strict.   
Intuitively, this condition guarantees that all the concentrations of species of a 
semilocking set can not simultaneously decrease  
while preserving the necessary monomial dominance.  D. Angeli {\em et al.}  
proved that if every semilocking set is dynamically non-emptiable, and 
if another condition holds (see \cite{Sontag07} for details), then the 
system is persistent. 
 
We next prove that equation \eqref{eq:bigbound} and a slight variant of 
the surrounding argument can be used to show that any semilocking set 
$W$ associated with a facet of a weakly reversible system is 
dynamically non-emptiable.  We therefore have provided a large set 
of examples of dynamically non-emptiable semilocking sets. 
 
\begin{cor} 
  Let $\{\S, \C, \Re\}$ be a weakly reversible chemical reaction 
  network with dynamics governed by mass action kinetics 
  \eqref{eq:main}.  Then any semilocking set associated with a facet 
  is dynamically non-emptiable. 
  \label{cor:dyn_non} 
\end{cor} 
 
\begin{proof} 
  As in the proof of Theorem \ref{thm:main}, we may assume that there is one linkage class.  The case of more than one linkage class is similar.  Let $W = \{X_1,\dots,X_M\}$ be a semilocking set such that $F_W$ is a facet.  Let $0 \ne \alpha \in \mathcal{F}_{\epsilon}(W)$ for some $\epsilon > 0$ which may be made smaller as needed.  Let $v \in \R^N_{\ge 0}$, with $v_i > 0$ if $i \in [M]$, be as in the proof of Theorem \ref{thm:main}; that is, $y'| _{W} - y|_{W} ~\in ~ \text{span}(v|_{W})$ for each reaction $y \to y' \in \Re$ and ${y'_k}|_W - {y_k}|_W = \gamma_k v|_{W}$ for some $\gamma_k \in \R$ where $k \in [R]$.  Then, for $w$ as in the definition of $C(W)$, and $R_{+}$ and $R_{-}$ defined similarly as in the proof of Theorem \ref{thm:main}, we have 
  \begin{align*} 
    w|_{W} ~=~ \sum_{k} \alpha_k(y_k'|_{W} - y_k|_{W}) ~=~ v|_{W}\sum_{y_k \to y_k' \in R_{+}} \alpha_k \gamma_k ~-~ v|_{W}\sum_{y_k \to y_k' \in R_{-}} \alpha_k |\gamma_k|~, 
  \end{align*} 
  and so for $i \in \{1,\dots,M\},$ 
  \begin{align} 
    w_{i} ~ = ~ v_{i}\sum_{y_k \to y_k ' \in R_{+}} \alpha_k \gamma_k ~-~ v_{i}\sum_{y_k \to y_k' \in R_{-}} \alpha_k |\gamma_k|~. 
    \label{eq:pos_neg} 
  \end{align} 
  Just as in the proof of Theorem \ref{thm:main}, there exists a reaction, the $\ell$th say,  $\tilde y_{\ell} \to y_{\ell}' \in R_{+}$ such that $\tilde{y}_{\ell} \preceq y$ for all $y \in \C$, and $\tilde{y}_{\ell} \ne y_k$ if $y_k \to y_k' \in R_{-}$.  Thus, for $\epsilon > 0$ small enough, and by the definition of $\mathcal{F}_{\epsilon}(W)$, the absolute value of the entire negative term in \eqref{eq:pos_neg} is less than or equal to $ v_{\text{min}} \alpha_{\ell} \gamma_{\ell}$, where $v_{\text{min}} := \min\{v_i \, : \, i \in [M]\} $.  Therefore we see that $w|_{W} \succeq 0$.  If $\alpha_{\ell} = 0$, then the preceding argument shows that $\alpha_k = 0$ for $y_k \to y'_k \in R_{-}$ and so $w_i > 0$ for all $i \in [M]$ since $\alpha \ne 0$; hence $\alpha \notin C(W)$.  For the case $\alpha_{\ell} \ne 0$, it is clear that $w_i > 0$ for all $i\in [M]$.  Thus, we again obtain that $\alpha \notin C(W)$ and the result is shown. 
  \end{proof}

\section{The Global Attractor Conjecture} 
\label{sec:GAC} 
 
In this section, we use the results of the previous section to resolve 
some special cases of the Global Attractor Conjecture of chemical 
reaction network theory.  In particular, the main result of this 
section, Theorem \ref{thm:facetVtx}, establishes that the conjecture 
holds if all boundary equilibria are confined to the facets and 
vertices of the positive stoichiometric compatibility classes.  That 
is, the conjecture holds if the faces associated with semilocking sets 
are facets, vertices, or are empty. 
 
\subsection{Complex-balancing systems} 
 
The Global Attractor Conjecture is concerned with the asymptotic 
stability of so-called ``complex-balancing'' equilibria.  Recall that 
a concentration vector $\overline x \in \R^N_{\geq 0}$ is an 
{\em equilibrium} of \eqref{eq:main} if the differential equations 
vanish at $\overline x$: $f(\overline x) = 0$. For each complex $\eta 
\in \C$ we will write $\{k \ | \ y_k = \eta\}$ and $\{k \ | \ y_k' = 
\eta\}$ for the subsets of reactions $k \in \mathcal{R}$ for which 
$\eta$ is the source and product complex, respectively.  In the 
following definition, it is understood that when summing over the 
reactions $\{k \ | \ y_k' = \eta\}$, $y_k$ is used to represent the 
source complex of the given reaction.

% :Definition of c-balancing 
\begin{defn} 
  An equilibrium $\overline x \in \R^N_{\ge 0}$ of \eqref{eq:main} is 
  said to be {\em complex-balancing} if the following equality  
  holds for each complex $\eta \in \C$: 
  \begin{equation*} 
    \sum_{\{k \ | \ y_k = \eta \}} \kappa_k   ( {\overline x} )^{y_k'} 
    ~=~ \sum_{\{k \ | \ y_k' = \eta\}}  \kappa_k ( {\overline x} ) 
    ^{y_k}  ~.  
  \end{equation*} 
  That is, $\overline x$ is a complex-balancing equilibrium if, at 
  concentration $\overline x$, the 
  sum of reaction rates for reactions for which  $\eta$ is the source is equal to 
  the sum of reaction rates for reactions for which  $\eta$ is the product. 
  A {\em complex-balancing system} is a mass action system 
  (\ref{eq:main}) that admits a strictly positive complex-balancing 
  equilibrium. 
  % need the assumption of strictly positive, because the presence of 
  % boundary c-bal equilibria does not imply the same for the interior 
\end{defn} 
 
In \cite{TDS}, complex-balancing systems are called ``toric dynamical 
systems'' in order to highlight their inherent algebraic structure. 
 Complex-balancing systems 
are automatically weakly reversible \cite{Feinberg72}.  There are two 
important special cases of complex-balancing systems: the 
detailed-balancing systems and the zero deficiency systems. 
 
% :Definition of d-balancing 
\begin{defn} 
  An equilibrium $\overline x \in \R^N_{\ge 0}$ of a reversible system 
  with dynamics given by mass action \eqref{eq:main} is said to be 
  {\em detailed-balancing} if for any pair of reversible reactions 
  $y_k \rightleftarrows y_k'$ with forward reaction rate $\kappa_k$ 
  and backward rate $\kappa_k'$, the following equality holds: 
  \begin{align*} 
    \kappa_k ( {\overline x} )^{y_k} = \kappa_k' ( {\overline x} 
    )^{y_k'}~. 
  \end{align*} 
  That is, $\overline x$ is a detailed-balancing equilibrium if the 
  forward rate of each reaction equals the reverse rate at 
  concentration $\overline x$.  A {\em detailed-balancing system} is a 
  reversible system with dynamics given by mass action \eqref{eq:main} that admits 
  a strictly positive detailed-balancing equilibrium. 
\end{defn} 
 
Properties of detailed-balancing systems are described by Feinberg in \cite{Fein89} and by   
A.~Vol\'{}pert and S.~Khud\t{ia}ev in \cite[Section 12.3.3]{Volp}. 
It is clear that detailed-balancing implies complex-balancing. 
 
\begin{defn} 
  For a chemical reaction network 
  $\{\mathcal{S},\mathcal{C},\mathcal{R}\}$, let $n$ denote the number 
  of complexes, $l$ the number of linkage classes, and $s$ the 
  dimension of the stoichiometric subspace, $S$.  The 
  {\em deficiency} of the reaction network is the integer $n-l-s$. 
\end{defn} 
 
The deficiency of a reaction network is non-negative 
because it can be interpreted as either the dimension of a 
certain linear subspace \cite{Fein79} or the codimension of a certain 
ideal \cite{TDS}.  Note that the deficiency depends only on the 
reaction network or the reaction diagram. It is known that any weakly 
reversible dynamical system \eqref{eq:main} whose deficiency is zero 
is complex-balancing, and that this fact is independent of the choice 
of rate constants $\kappa_k$ 
\cite{Fein79}.   
% In fact, this property of being c-bal regardless of rate parameters 
% defines the space of deficiency zero systems. 
On the other hand, a reaction diagram with a deficiency that is 
positive may give rise to a system that is both complex- and 
detailed-balancing, complex- but not detailed-balancing, or neither, 
depending on the values of the rate constants $\kappa_{k}$ 
\cite{TDS,Feinberg72,Fein89,Horn72}. 
 
\subsection{Qualitative behavior of complex-balancing systems}  
\label{sec:qual} 
%and 
%the Global Attractor Conjecture} 
 
Much is known about the limiting behavior of complex-balancing 
systems.  In the interior of each positive stoichiometric 
compatibility class $P$ for such a system, there exists a unique equilibrium 
$\overline{x}$, with strictly positive components, and this 
equilibrium is complex-balancing.  As previously stated, $\overline x$ 
is called the {\em Birch point} due to the connection to Birch's 
Theorem (see Theorem~1.10 of \cite{ASCB}).  Note that a system was 
defined to be complex-balancing if at least one such equilibrium 
exists; we now are asserting that so long as at least one $P$ contains 
a complex-balancing equilibrium, then they all do.  As for the 
stability of the equilibrium within the interior of the corresponding 
$P$, a strict Lyapunov function exists for each such point.  Hence 
local asymptotic stability relative to $P$ is guaranteed; see Theorem 
6A of \cite{HornJackson} and the Deficiency Zero Theorem~of 
\cite{Fein79}.  %Note - Lyapunov function not infinity on boundary. 
The Global Attractor Conjecture %(GAC) 
states that this equilibrium of $P$ is globally asymptotically stable 
relative to the interior of $P$ \cite{TDS}.  In the following 
statement, a {\em global attractor} for a set $V$ is a point $v^{*} 
\in V$ such that any trajectory $v(t)$ with initial condition $v^{0} 
\in V$ converges to $v^{*}$, in other words, $\lim\limits_{t 
  \rightarrow \infty} v(t) = 
v^{*}$. %check H & Smale for precise definition? 
 
% :GAC statement 
\medskip 
 
\noindent \textbf{Global Attractor Conjecture} {\em For any 
  complex-balancing system \eqref{eq:main} and any strictly positive 
  initial condition $x^0$, the Birch point $\overline x \in P := (x^0 
  + S) \cap \R_{\ge 0}^N$ is a global attractor of the interior of the 
  positive stoichiometric compatibility class, $\operatorname{int} (P)$.  } 
 
\medskip This conjecture first appeared in a paper of F. Horn 
\cite{Horn74}, and was given the name ``Global Attractor Conjecture'' 
by G. Craciun {\em et al.} \cite{TDS}.  It is stated to be the main open 
question in the area of chemical reaction network theory by L. Adleman 
{\em et al.}  \cite{Adleman08}.  In fact, M. Feinberg stated the more  
general conjecture that all weakly reversible systems are persistent  
\cite[Section 6.1]{Fein79}.  To this end, G. Gnacadja proved 
that the class of networks of ``reversible binding reactions'' are 
persistent; these systems include non-complex-balancing 
ones~\cite{Gnacadja08}. 
 
We now describe known partial results regarding the \GACnospace.  By 
an {\em interior trajectory} we shall mean a solution $x(t)$ to 
(\ref{eq:main}) that begins at a strictly positive initial condition 
$x^0 \in \mathbb{R}^N_{>0}$.  It is known that trajectories of complex-balancing systems converge to the set of equilibria; see Corollary 
2.6.4 of \cite{ChavezThesis} or Theorem~1 of \cite{Sontag01}. 
% so only one equilibrium to converge to. 
Hence, the conjecture is equivalent to the following statement: {\em 
  for a complex-balancing system, any boundary equilibrium is not an 
  $\omega$-limit point of an interior trajectory}.  It clearly follows 
that if a positive stoichiometric compatibility class $P$ has no 
boundary equilibria, then the \GAC holds for this $P$.  %This is Theorem 4.1 
 % of Siegel and McLean 
Thus, sufficient conditions for the non-existence of boundary 
equilibria are conditions for which the \GAC holds (see Theorem~2.9 of 
\cite{Anderson08}); a result of this type is Theorem~6.1 of L.~Adleman 
{\em et al.}  \cite{Adleman08}.  Recall that by 
Theorem~\ref{thm:omegapoints}, we know that the only faces $F_{W}$ of 
a positive stoichiometric compatibility class $P$ that may contain 
$\omega$-limit points in their interiors are those for which $W$ is a 
semilocking set.  In particular, if the set $Z_{W}$ is 
stoichiometrically unattainable for all semilocking sets $W$, then $P$ 
has no boundary equilibria, and hence, the \GAC holds for this $P$; 
see the main theorem of D. Angeli {\em et al.}  \cite{PetriNet}. 
Biological models in which the non-existence of boundary equilibria 
implies global convergence include the ligand-receptor-antagonist-trap 
model of G. Gnacadja {\em et al.} \cite{LRAT}, the enzymatic mechanism 
of D. Siegel and D. MacLean \cite{SiegelMacLean}, and the T-cell 
signal transduction model of T.~McKeithan \cite{T-cell} (the mathematical analysis 
appears in the work of E. Sontag \cite{Sontag01} and in 
the Ph.D. thesis of M.  Chavez \cite[Section~7.1]{ChavezThesis}).  We remark that 
this type of argument first appeared in~\cite[Section~6.1]{FeinDefZeroOne}. 
%see also \cite{Anderson08, PetriNet}. 
% Maybe also cite the two Siegel ``Global stability'' papers?? 
 
The remaining case of the \GACnospace, in which equilibria exist on 
the boundary of $P$, is still open.  However, some progress has been 
made.  For example, it already is known that vertices of $P$ can not 
be $\omega$-limit points even if they are equilibria; see Theorem~3.7 
in the work of D. Anderson \cite{Anderson08} or Proposition~20 of the 
work of G. Craciun {\em et al.} \cite{TDS}.  For another class of 
systems for which the \GAC holds despite the presence of boundary 
equilibria, see Proposition 7.2.1 of the work of M. Chavez 
\cite{ChavezThesis}.  The hypotheses of this result are that the set 
of boundary equilibria in $P$ is discrete, that each boundary 
equilibrium is hyperbolic with respect to $P$, and that a third, more 
technical condition holds. 
%(H2) holds [must read what this is!] 
In addition, the \GAC holds in the case that the network is 
detailed-balancing, $P$ is two-dimensional, and the network is {\em 
  conservative} (meaning that $P$ is bounded); see Theorem~23 of 
\cite{TDS}.  In the next section, Corollary \ref{cor:2d} will allow us to 
eliminate the hypotheses ``detailed-balancing'' and ``conservative'' 
from the two-dimensional result; such partial results concerning the \GAC  
are the next topic of this paper. 
 
\subsection{Applications to complex-balancing systems}

Theorem \ref{thm:facetVtx} is our main contribution to the Global 
Attractor Conjecture and generalizes the known results described in 
Section \ref{sec:qual}.  We first present a definition and a lemma. 
 
\begin{defn} 
  Suppose that $\{\S,\C ,\Re\}$ is a weakly reversible chemical 
  reaction network, endowed with mass action kinetics, and $W \subset 
  \S$ is a semilocking set.  Then, the {\em $W$-reduced network} is 
  the chemical reaction network $\{\S \setminus W,\tilde{\C}, \tilde{\Re}\}$ such 
  that $\tilde{\C}$ and $\tilde{\Re}$ are those complexes   
  and 
  reactions that do not involve a species from $W$.  Equivalently, this is the subnetwork obtained by removing all linkage classes that contain a complex 
  with a species from $W$ in its support.  The {\em $W$-reduced 
    system} is the $W$-reduced network endowed with the same rate 
  constants as in the original system. 
\end{defn} 
 
As noted in comments following Proposition~\ref{prop:semilocking_intution}, for a 
weakly reversible system and any semilocking set $W$, either each 
complex in a given linkage class contains an element of $W$ or each 
complex in that linkage class does not contain any element of $W$. 
Therefore, $W$-reduced systems % are well defined, and 
are themselves weakly reversible.  % b/c {S,C,R} is weakly reversible 
 Furthermore, it is easy to check that any $W$-reduced system of a complex-balancing system is itself complex-balancing. 
 
\begin{lem} 
  Consider a complex-balancing system.  A face $F_{W}$ of a 
  stoichiometric compatibility class $P$ contains an equilibrium in 
  its interior if and only if $W$ is a semilocking set. 
  \label{lem:classification} 
\end{lem} 
 
\begin{proof} 
  It is clear that if a face $F_W$ contains an equilibrium in its 
  interior, then $W$ is a semilocking set (see the discussion 
  following Proposition~\ref{prop:semilocking_intution} or the proof of Theorem 2.5 
  of \cite{Anderson08}). 
 
  If $W$ is a semilocking set, then the $W$-reduced system is 
  complex-balancing and one of its invariant polyhedra is naturally identified with $F_W$.    
  Therefore, this reduced system admits a complex-balancing 
  equilibrium, which we denote by $y$.  For $i \in W^c$, let $z_i$ be the component of  
  $y$ associated with species $i$.  For $i \in W$, let $z_i = 0$.  So 
  constructed, the concentration vector $z$ is an equilibrium within the interior of $F_W$. 
\end{proof}

\begin{thm} 
  The Global Attractor Conjecture holds for any complex-balancing (and 
  in particular, detailed-balancing or weakly reversible zero 
  deficiency) chemical reaction system whose boundary equilibria are 
  confined to facet-interior points or vertices of the positive 
  stoichiometric compatibility classes.  Equivalently, if a face 
  $F_{W}$ is a facet, vertex, or an empty face whenever $W$ is a 
  semilocking set, then the \GAC holds. 
  \label{thm:facetVtx} 
\end{thm} 
 
\begin{proof} 
  The equivalence of the two statements in the theorem is a consequence of  
  Lemma~\ref{lem:classification}.  As noted in the previous section, 
  persistence is a necessary and sufficient condition for the Global 
  Attractor Conjecture to hold.  Further, by the results in 
  \cite{Anderson08} or \cite{TDS}, vertices may not be $\omega$-limit 
  points.  The remainder of the proof is similar to that of 
  Theorem~\ref{thm:persFacet} and is omitted. 
\end{proof} 
 
We note that Theorem~\ref{thm:facetVtx} and previous results in \cite{Anderson08,TDS} show that for equilibria $z$  
of complex-balancing systems that reside within the interiors of 
facets or of vertices of $P$, these faces are repelling near $z$.  More precisely, there exists a relatively open set $Q$ of $P$ such that $Q$ contains $z$ and the face is repelling in $Q\cap \operatorname{int}(P)$. 
%For boundary equilibria on faces of dimensions in between the 
%two extreme cases, repelling neighborhoods may not exist.  An example 
%of this phenomenon appears in Example \ref{ex:4}. 
 
The following corollary % to Theorem~\ref{thm:facetVtx} 
resolves the \GAC for systems of dimension two; note that the 
one-dimensional case is straightforward to prove. 
 
\begin{cor}[GAC for two-dimensional $P$]\label{cor:2d} 
  The Global Attractor Conjecture holds for all complex-balancing (and 
  in particular, detailed-balancing or weakly reversible zero 
  deficiency) chemical reaction systems whose positive stoichiometric 
  compatibility classes are two-dimensional. 
\end{cor} 
 
\begin{proof} 
  This follows immediately from Theorem~\ref{thm:facetVtx} as each 
  face of a two-dimensional polytope must be either a facet or a 
  vertex. 
 \end{proof} 
 
 In the next section, we provide examples that illustrate our results, 
 as well as a three-dimensional example for which our results do not 
 apply. 
%for which the \GAC is still open. 

% :Section: EXAMPLES 
\section{Examples}\label{sec:exs} 
As discussed in the previous section, the \GAC previously has been 
shown to hold if $P$ has no boundary equilibria or if the boundary 
equilibria are restricted to vertices of $P$.  Therefore, the examples 
in this section that pertain to complex-balancing systems feature 
non-vertex boundary equilibria.  Also, we have chosen  
examples for which the conditions of the theorems can be easily checked. 
 
\begin{example}\label{ex:1} 
  We revisit the network given by the following reactions: 
  \begin{align*} 
    2A ~ \underset{\kappa_2}{\overset{\kappa_1}{\rightleftarrows}} ~ 
    A+B \quad , \quad B ~ 
    \underset{\kappa_4}{\overset{\kappa_3}{\rightleftarrows}} ~ C ~. 
  \end{align*} 
  As we saw in Example \ref{ex:P}, the positive stoichiometric 
  compatibility classes are two-dimensional triangles: 
  \begin{align} 
    P~=~ \left\{~ (x_a, x_b, x_c ) \in \mathbb{R}^3_{\geq 0} ~|~ 
      x_a+x_b + x_c = T ~ \right\}~, 
    \label{triangle} 
  \end{align} 
  where  
  $T>0$.   
  % The three facets (edges) of each positive compatibility class $P$ 
  % are $F_{\{A\}}$, $F_{\{B\}}$, and $F_{\{C\}}$. 
  It is straightforward to check that $P$ has a unique boundary 
  equilibrium given by 
  \begin{align*} 
    z~=~ \left( 0, ~ \frac{\kappa_4 }{\kappa_3 +\kappa_4} T, ~ 
      \frac{\kappa_3 }{\kappa_3 +\kappa_4} T \right)~, 
  \end{align*} 
  and that this point lies in the interior of the facet $F_{ \{A\} } 
  $.  (Note that this boundary equilibrium is the Birch point of the 
  reversible deficiency zero subnetwork $B \leftrightarrows C$.) 
  Therefore, both Theorem \ref{thm:facetVtx} and Corollary 
  \ref{cor:2d} allow us to conclude that despite the presence of the 
  boundary equilibrium $z$, the Birch point in the interior of $P$ is 
  globally asymptotically stable. 
  %[Check: does M. Chavez thesis result apply?] 
\end{example} 
 
We remark that the results in \cite{Anderson08} do not apply to the 
previous example, although Theorem~23 of \cite{TDS} and Theorem~4 of 
\cite{Sontag07} do.  However, for the following example, no previously 
known results apply. 
 
% :Example2 
%\begin{figure} 
%  \begin{center} 
%    \includegraphics[scale=0.9]{newfigure.eps} 
%  \end{center} 
%  \caption{Chemical reaction network for Example \ref{ex:2}.} 
%   \label{triangleNetwork} 
%\end{figure} 
 
\begin{example}\label{ex:2} 
  Consider the reaction network depicted below 
  \begin{center} 
  \includegraphics[scale=0.9]{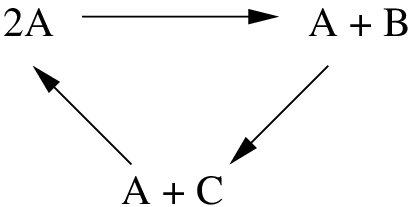}. 
  \end{center} 
  The positive stoichiometric compatibility classes are the same 
  triangles (\ref{triangle}) as in the previous example.  For each 
  $P$, the set of boundary equilibria is the entire face $F_{ \{ {A} 
    \} }$ (one of the three edges of $P$), which includes the two 
  vertices $F_{\{A,B\}}$ and $F_{\{A,C\}}$.  Hence the results of 
  \cite{Anderson08,ChavezThesis} do not apply.  Note that this is a 
  weakly reversible zero deficiency network, but it is not 
  detailed-balancing; so the results of \cite{TDS} do not apply. 
  Finally, the second condition of Theorem 4 of \cite{Sontag07} is not 
  satisfied by this example (there are nested ``deadlocks'' 
  \cite{Sontag07}) and so that result also does not apply. 
  However, both Theorem \ref{thm:facetVtx} and Corollary \ref{cor:2d} 
  imply that the \GAC holds for all choices of rate constants and for 
  all $P$ defined by this network, despite the presence of boundary 
  equilibria. 
\end{example} 
 
In the next example, the positive stoichiometric compatibility classes 
are three-dimensional. 
% :Example 3 
\begin{example}\label{ex:3} 
  The following zero deficiency network is obtained from 
  Example~\ref{ex:1} by adding a reversible reaction: 
  \begin{align*} 
    2A ~ \leftrightarrows ~ A+B~ \quad ,  \quad \quad B ~ 
    \leftrightarrows ~C~ \leftrightarrows ~ D~. 
  \end{align*} 
  The positive stoichiometric compatibility classes are 
  three-dimensional simplices (tetrahedra): 
  \begin{align*} 
    P~=~ \left\{~ (x_a, x_b, x_c,x_{d} ) \in \mathbb{R}^4_{\geq 0} ~|~ 
      x_a+x_b + x_c +x_{d} = T ~ \right\}~, 
  \end{align*} 
  for positive total concentration $T>0$.  The unique boundary 
  equilibrium in $P$ is the Birch point of the zero deficiency 
  subnetwork $B \leftrightarrows C \leftrightarrows D$, and it lies in 
  the interior of the  facet $F_{\{A \} }$.  In other words, the point is $z=(0,x_b, 
  x_c,x_{d})$ where $(x_b, x_c,x_{d})$ is the Birch point for the 
  system defined by the subnetwork 
  \begin{equation*} 
  B ~   \leftrightarrows ~C~ \leftrightarrows ~ D~. 
  \end{equation*} 
 Thus by Theorem~\ref{thm:facetVtx} %\ref{thm:facetVtx}, 
  the \GAC holds for all positive stoichiometric compatibility classes $P$ and all choices of rate constants defined 
  by this network. 
  % [Check: does Proposition 7.2.1. of M. Chavez's thesis apply 
  % \cite{ChavezThesis}?] 
\end{example} 
 
As in the previous example, the positive stoichiometric compatibility 
classes of our next example are three-dimensional.  However neither 
previously known results \cite{Anderson08,ChavezThesis,TDS} nor our 
current results can resolve the question of global asymptotic 
stability. 
% :Example 4 
\begin{example}\label{ex:4} 
  The following zero deficiency network consists of three reversible 
  reactions: 
  \begin{align*} 
    A ~ \leftrightarrows ~ B ~ \leftrightarrows ~ A+B ~ 
    \leftrightarrows ~A+C~. 
  \end{align*} 
  %haven't defined ``conservation relations,'' but that's fine for now 
  As there are no conservation relations, the unique positive 
  stoichiometric compatibility class is the entire non-negative 
  orthant: 
  \begin{align*} 
    P~=~ \mathbb{R}^3_{\geq 0} ~. 
  \end{align*} 
  The set of boundary equilibria is the one-dimensional face (ray) 
  $F_{\{A,B\}}$, which includes the origin $F_{\{A,B,C\}}$.  Therefore 
  non-vertex, non-facet boundary equilibria exist, so the results in 
  this paper do not apply. 
  \end{example} 
   
  We end with an example in which the results of 
  Section~\ref{sec:MainThm} apply but those of Section~\ref{sec:GAC} 
  do not. 
   
%:Example 5 
  \begin{example} 
    The following reversible network is obtained from Example 
    \ref{ex:1} by adding another reversible reaction: 
  \begin{align*} 
    A + C ~\underset{\kappa_2}{\overset{\kappa_1} {\leftrightarrows}} 
    ~ 2A ~ \underset{\kappa_4}{\overset{\kappa_3}{\rightleftarrows}} ~ 
    A+B \quad , \quad B ~ 
    \underset{\kappa_6}{\overset{\kappa_5}{\rightleftarrows}} ~ C ~. 
  \end{align*} 
  The positive stoichiometric compatibility classes are again the 
  two-dimensional triangles given by \eqref{triangle}.  
%  \begin{align*} 
    %P~=~ \left\{~ (x_a, x_b, x_c ) \in \mathbb{R}^3_{\geq 0} ~|~ 
      %x_a+x_b + x_c = T ~ \right\}~, 
  %\end{align*} 
  %where $T>0$.   
One can easily check that the network has a deficiency 
  of one, so there exist rate constants for which the system is not complex-balancing (for example, $\kappa_1 = \kappa_3$, 
  $\kappa_5 = \kappa_6$, and $\kappa_2 \ne \kappa_4$).  Thus the 
  results of Section \ref{sec:GAC} do not apply.  It is also easy to 
  verify that $\{A\}$ and $\{A,B,C\}$ are the only semilocking sets and 
  that $F_{\{A\}}$ is a facet and $F_{\{A,B,C\}}$ is empty. 
  Therefore, Theorem \ref{thm:persFacet} applies and we conclude that, 
  independent of the choice of rate constants, the system is 
  persistent. 
    \end{example} 
 
% :Acknowledgments 
    \subsubsection*{Acknowledgments} The authors thank Gheorghe 
    Craciun for helpful discussions.  This work began during a focused 
    research group hosted by the Mathematical Biosciences Institute at 
    The Ohio State University, and the authors benefited from a 
    subsequent visit to the Statistical and Applied Mathematical 
    Sciences Institute in North Carolina.   
    We also acknowledge the helpful comments of anonymous reviewers, which greatly improved the paper. 
 
% :Bibliography 
% \bibliographystyle{amsplain} \bibliography{Facet} 

\end{document}